\documentclass [11pt]{amsart}
\usepackage[ansinew]{inputenc}
\usepackage{amssymb}
\usepackage{amsmath}
\usepackage{graphicx}
\usepackage{amsthm}
\usepackage{cite}
\usepackage{color}

\newtheorem{theorem}{\textbf{Theorem}}[section]
\newtheorem{lemma}{\textbf{Lemma}}[section]
\newtheorem{proposition}{\textbf{Proposition}}[section]
\newtheorem{corollary}{\textbf{Corollary}}[section]
\newtheorem{remark}{\textbf{Remark}}[section]
\newtheorem{definition}{\textbf{Definition}}[section]

\allowdisplaybreaks[4]

\def\be{\begin{equation}}
\def\ee{\end{equation}}
\def\bea{\begin{eqnarray}}
\def\eea{\end{eqnarray}}
\def\bt{\begin{theorem}}
\def\et{\end{theorem}}
\def\bl{\begin{lemma}}
\def\el{\end{lemma}}
\def\br{\begin{remark}}
\def\er{\end{remark}}
\def\bp{\begin{proposition}}
\def\ep{\end{proposition}}
\def\bc{\begin{corollary}}
\def\ec{\end{corollary}}
\def\bd{\begin{definition}}
\def\ed{\end{definition}}

\def\non{\nonumber }

\DeclareMathOperator{\tr}{tr}

\newcommand{\F}{\mathcal F}

\newcommand{\Id}{\mathbb I}

\newcommand\defeq{\stackrel{\scriptscriptstyle \text{def}}=}
\newcommand{\RR}{\mathbb R}

\newcommand{\sS}{\mathcal{S}}

\newcommand \dQ{Q^\delta}
\newcommand \dR{R^\delta}

\usepackage{color}
\definecolor{red}{rgb}{1,0,0}

\definecolor{blue}{rgb}{0,0,1}

\begin{document}

\title[eigenvalues for the co-rotational Beris-Edwards system]{An elementary proof of eigenvalue preservation for the co-rotational Beris-Edwards system}

\author{ Andres Contreras}
\thanks{
  Department of Mathematical Sciences,
  New Mexico State University, Las Cruces, New Mexico, 88003, USA.
  \texttt{acontre@nmsu.edu}}

 \author{ Xiang Xu}
\thanks{
  Department of Mathematics and Statistics,
  Old Dominion University, Norfolk, Virginia, 23529, USA.
  \texttt{x2xu@odu.edu}}

 \author{ Wujun Zhang}
\thanks{
Department of Mathematics, Rutgers University, Piscataway, New Jersey,
08854, USA.
  \texttt{wujun@math.rutger.edu}}

\date{}
\maketitle

\begin{abstract}
We study the co-rotational Beris-Edwards system modeling nematic
liquid crystals and revisit the eigenvalue preservation property
discussed in \cite{XZ16}. We give an alternative but direct proof to
the eigenvalue preservation of the initial data for the $Q$-tensor.
It is noted that our proof is not only valid in the whole space
case, but in the bounded domain case as well.

\end{abstract}

\section{Introduction}
\setcounter{equation}{0}

In this paper we study the eigenvalue preservation property of
solutions for a hydrodynamic system modeling the evolution of
nematic liquid crystals. Mathematically speaking, this system is
composed of a coupled incompressible Navier-Stokes equations with
anisotropic forces and Q-tensor equations of a parabolic type that
describes the evolution of the liquid crystal director field, which
is called the Beris-Edwards system \cite{BE94}. In the Landau-de
Gennes theory \cite{B12, dG93}, the basic element is a symmetric,
traceless tensor $Q$ that is a tensor valued function taking values
in the $5$-dimensional $Q$-tensor space
$$
  \sS_0^{(3)}\defeq\big\{Q\in\mathbb{M}^{3\times 3}, \, Q^{t}=Q, \, \tr(Q)=0\big\}.
$$

The simplest form of the free energy in the Landau-de Gennes theory
takes the following form:
\begin{equation}\label{free energy}
\F(Q)\defeq
\int_{\Omega}\frac{L}{2}|\nabla{Q}|^2+\frac{a}{2}\tr(Q^2)-\frac{b}{3}\mathrm{tr}(Q^3)+\frac{c}{4}\tr^2(Q^2)\,dx,
\end{equation}
where $\Omega\subset\mathbb{R}^3$ is a smooth and bounded domain.
Above in \eqref{free energy} we use the one constant approximation
of the Oseen-Frank energy, and $L, a, b, c$ are material dependent
constants that satisfy \cite{M10, MZ10}
\begin{equation}\label{AssumptionsBulk}
L>0, \; b> 0, \; c>0.
\end{equation}

The simplified Beris-Edwards system we study reads
\begin{align}
\mathbf{u}_t+\mathbf{u}\cdot\nabla\mathbf{u}-\nu\Delta\mathbf{u}+\nabla{P}&=
\lambda{L}\nabla\cdot(Q\Delta{Q}-\Delta{Q}Q)-\lambda{L}\nabla\cdot(\nabla{Q}\odot\nabla{Q}),\label{NSE}\\
\nabla\cdot\mathbf{u}&=0,\label{incomp}\\
Q_t+\mathbf{u}\cdot\nabla{Q}-\omega{Q}+Q\omega&=\Gamma\Big(L\Delta{Q}-aQ+b\Big[Q^2-\frac{\tr(Q^2)}{3}\mathbb{I}\Big]-cQ\tr(Q^2)\Big),
\label{Q-equ}
\end{align}
with the following initial and boundary conditions
\begin{align}\label{IC-BC}
&\mathbf{u}(0, x)=\mathbf{u}_0(x) \;\mbox{with} \;
\nabla\cdot\mathbf{u}_0=0, \quad Q(0, x)=Q_0(x)\in \sS_0^{(3)},\non\\
&\mathbf{u}(t, x)|_{\partial\Omega}=0, \quad
Q(t,x)|_{\partial\Omega}=Q_0(x)|_{\partial\Omega}=\tilde{Q}(x).
\end{align}
Above $\mathbf{u}(t, x): (0,
+\infty)\times\Omega\rightarrow \RR^3$ stands for the incompressible
fluid velocity field, $Q(t, x): (0, +\infty)\times\Omega\rightarrow
\sS_0^{(3)}$ represents the order parameter of the liquid crystal
molecules and
$\omega=\dfrac{\nabla\mathbf{u}-\nabla^T\mathbf{u}}{2}$ denotes the
skew-symmetric part of the rate of strain tensor. The positive
constants $\nu, \lambda$ and $\Gamma$ denote the fluid viscosity,
the competition between kinetic energy and elastic potential energy,
and macroscopic elastic relaxation time for the molecular
orientation field, respectively \cite{XZ16}. This simplified system
is at time referred to as the ``co-rotational" Beris-Edwards system
\cite{PZ12} in the literature, whose related mathematical study can
be found in \cite{ADL15, D15, DFRSS14, GR14, GR15, PZ12}. On the
other hand, the full system is also called the ``non co-rotational"
Beris-Edwards system, and we refer interested readers to
\cite{ADL14, CRWX15, DAZ, PZ11, Z17, Z17-2} for its relevant PDE and
numeric work.

From the physical point of view, the main feature of nematic liquid
crystals is the locally preferred orientation of the nematic
molecule directors. To this end $Q$-tensors are introduced, which
are considered suitably normalized second order moments of the
probability distribution function. Specifically, if $\mu_x$ is a
probability measure on the unit sphere $\mathbb{S}^2$ representing
the orientation of liquid crystal  molecules at a point $x$ in
space, then a $Q$-tensor denoted by ${Q}(x)$ is a symmetric and
traceless $3\times 3$ matrix defined by
\begin{equation}\label{de Gennes Q-tensor}
  Q(x)=\int_{\mathbb{S}^2}\left(\mathbf{p}\otimes  \mathbf{p}-\frac{1}{3}\Id\right)\,d\mu_x(\mathbf{p}).
\end{equation}
Indeed it is a crude measure (from the viewpoint of statistical
theory) of how the second-moment tensor associated with a given
probability measure deviates from its isotropic value
\cite{NMreview, V94}. It is noted that \eqref{de Gennes Q-tensor}
imposes a constraint such that (see \cite{NMreview})
$$
-\frac{1}{3}\leq \lambda_i(Q)\leq  \frac{2}{3},\quad\forall\, 1\leq
i\leq 3.
$$
Hence it is easy to check that not every symmetric and traceless
$3\times 3$ matrix is a {\it physical} $Q$-tensor but only those
whose eigenvalues range in $[-\frac{1}{3},\frac{2}{3}]$.

Motivated by the physical interpretation of the Q-tensors, it seems
to be of great importance to understand how the fluid dynamics would
affect the behavior of eigenvalues of the $Q$-tensors as time
evolves. Partially motivated by this question, in \cite{XZ16}, the
authors proved that certain eigenvalue constraints of the initial
data $Q_0$ are preserved by the evolution problem
\eqref{NSE}-\eqref{IC-BC} when the domain is either the entire
Euclidean space or a periodic box. Inspired by the idea in
\cite{XZ16}, in this paper we give an alternative but direct proof
whose argument works well both in the whole space case and in the
bounded domain case.

Our main result is stated as follows.

\begin{theorem}\label{main theorem}
For any $\mathbf{u}_0\in H_0^1(\Omega)$,
$\nabla\cdot\mathbf{u}_0=0$, $Q_0\in H^2(\Omega; \sS_0^{(3)})$ and
$\tilde{Q}\in H^\frac52(\partial\Omega)$, let $(\mathbf{u}(t,x),
Q(t,x))$ be the unique local strong solution to the evolution
problem \eqref{NSE}-\eqref{IC-BC} on $[0, T]$. We assume
\begin{equation}\label{assumption-a}
  0\leq a\leq\frac{b^2}{24c},
\end{equation}
and the initial data $Q_0$ and the boundary data  $\tilde{Q}$
satisfy
\begin{equation}\label{initial setting}
\lambda_i(Q_0(x))\in\Big[-\frac{b+\sqrt{b^2-24ac}}{12c},\frac{b+\sqrt{b^2-24ac}}{6c}\Big],
\qquad \forall x\in\Omega,\; 1\leq i\leq 3.
\end{equation}
Then for any $t\in (0,T]$ and
$x\in\Omega$, the eigenvalues of $Q(t,x)$ stay in the same interval.
\end{theorem}
\begin{remark}\label{remark on regularity}
By Theorem 1.1 in \cite{LY16}, the existence and uniqueness of
local strong solutions to the evolution problem
\eqref{NSE}-\eqref{IC-BC} is ensured, and satisfies
\begin{align*}
&\mathbf{u}\in H^1(0, T; H^1(\Omega))\cap L^\infty(0, T;
H^2(\Omega)), \quad\nabla\cdot\mathbf{u}=0, \\
&Q\in H^2\big(0, T; L^2(\Omega; \sS_0^{(3)})\big)\cap H^1\big(0, T;
H^2(\Omega; \sS_0^{(3)})\big)\cap L^\infty\big(0, T; H^3(\Omega;
\sS_0^{(3)})\big).
\end{align*}
\end{remark}

\begin{remark}
Compared to \cite{XZ16}, one extra assumption in Theorem \ref{main
theorem} is $a\geq 0$ which captures a regime of physical interest
but not the deep nematic regime \cite{DZ}. This assumption is only
used to get the same lower bound $-\frac{b+\sqrt{b^2-24ac}}{12c}$,
but not needed to achieve the upper bound
$\frac{b+\sqrt{b^2-24ac}}{6c}$. We also want to point out that this
assumption \eqref{assumption-a} is different from its counterpart
imposed in \cite{XZ16} because the bulk part are dealt with in
different ways.
\end{remark}

The idea of the proof is to proceed by contradiction and to exploit
the variational characterization of the eigenvalues in relation to
the evolution problem \eqref{NSE}-- \eqref{IC-BC}, which works for
the solution $Q$ with $C^{1, 2}$ regularity. If we were able to show
the solutions to the Beris-Edwards system were regular enough we
would be done, however this seems out of reach at the moment, though
an interesting problem on its own. Fortunately, we can bypass this
difficulty by using a regularization argument discussed in
\cite{XZ16} that preserves the eigenvalue constraints (the
eigenvalues converge pointwise in fact in the whole domain).

{For simplicity we set the eigenvalues of matrix $Q$
$$
  \lambda_i(t, x) \defeq\lambda_i(Q(t, x)), \quad 1\leq i\leq 3.
$$
Without loss of generality we assume
$$
  \lambda_1(t, x)\geq\lambda_2(t, x)\geq\lambda_3(t, x),
  \quad\forall (t, x)\in\bar\Omega\times [0, T]
$$}
{As a matter of fact, we may
establish the following more general result based on Theorem
\ref{main theorem}.}
\begin{corollary}\label{Corollary}
{For any given $\mathbf{u}_0\in H_0^1(\Omega)$,
$\nabla\cdot\mathbf{u}_0=0$, $Q_0\in H^2(\Omega; \sS_0^{(3)})$ and
$\tilde{Q}\in H^\frac52(\partial\Omega)$, the unique local strong
solution $(\mathbf{u}(t,x), Q(t,x))$ to the evolution problem
\eqref{NSE}-\eqref{IC-BC} on $[0, T]$ satisfies
\begin{align}
\lambda_1(t,x)&\leq\max\Big[\frac{b+\sqrt{b^2-24ac}}{6c},
\displaystyle\max_{\bar{\Omega}}\lambda_1(Q_0) \Big],\\
\lambda_3(t,x)&\geq\min\Big[-\frac{b+\sqrt{b^2-24ac}}{12c},
\displaystyle\min_{\bar{\Omega}}\lambda_3(Q_0) \Big],
\end{align}
for any $t\in (0,T]$ and $x\in\Omega$.}

\end{corollary}

\begin{remark} {Analogously,
Theorem \ref{main theorem} and Corollary \ref{Corollary} also valid
in the static case, provided the corresponding solution $Q\in
C(\bar{\Omega})\cap C^2(\Omega)$. However, this regularity issue
cannot be solved directly by following the approximation argument in
the appendix part, and henceforth is beyond the scope of our paper.}

\end{remark}

The proof of Theorem \ref{main theorem} {and Corollary
\ref{Corollary}} is given in Section 2, while a related technical
regularization lemma is presented in the appendix.

%


\section{Proof of Theorem \ref{main theorem}}
\setcounter{equation}{0} Here and after, we let $|Q|$ denote the
Frobenius norm of $Q\in\sS_0^{(3)}$, that is, $|Q| = \sqrt{tr(Q^t
Q)}$ where $tr$ denotes the trace of a matrix. Also, because of the
traceless property of $Q$-tensors, for all $(t,
x)\in\bar\Omega\times [0, T]$ one has

\begin{equation}\label{eigsumzero}
 \lambda_1(t, x)+\lambda_2(t, x)+\lambda_3(t, x)=0.
\end{equation}
To begin with, we see from \cite{R69, N73} that
\begin{lemma}\label{continuous eigenvalue}
For $1\leq i\leq 3$, $\lambda_i(t, x)\in C(\bar\Omega\times[0, T])$
\end{lemma}

Now we are ready to prove Theorem \ref{main theorem}.
\begin{proof}
Due to Lemma \ref{lemma on regularizing effect} in the Appendix, we
may assume
$$
Q(t, x)\in C^{1,2}((0,T)\times\Omega)\cap C([0, T]\times\bar\Omega).
$$

\noindent{\bf Step 1.} Let
\begin{equation}\label{maximum point}
\lambda_1(t_0, x_0)=\displaystyle\max_{(t, x)\in[0,
T]\times{\bar\Omega}}\lambda_1(t, x).
\end{equation}
We shall show that $ \lambda_1(t_0, x_0) \leq
\dfrac{b+\sqrt{b^2-24ac}}{6c}. $ We prove by a contradiction
argument. Suppose
\begin{equation}\label{large assumption}
(t_0, x_0)\in(0, T]\times\Omega, \;\mbox{and}\;\lambda_1(t_0,
x_0)>\frac{b+\sqrt{b^2-24ac}}{6c}.
\end{equation}
Let $\vec{v}\in\mathbb{S}^2$ be the corresponding unit eigenvector,
such that $Q(t_0,x_0)\vec{v}=\lambda_1(t_0,x_0)\vec{v}$. Meanwhile,
we denote
$$
  f(t, x)=\langle Q(t, x)\vec{v}, \vec{v}\rangle_{\mathbb{R}^3},
$$
then it is easy to check from \eqref{initial setting} and
\eqref{maximum point} that
\begin{equation}\label{maximum f}
f(t_0, x_0)=\displaystyle\max_{(t, x)\in[0,
T]\times{\bar\Omega}}f(t, x)
\end{equation}
Next, we take the matrix inner product of equation \eqref{Q-equ}
with $\vec{v}\vec{v}^t$ and evaluate the resultant at $(t_0,
x_0)$. Note that $\mathbf{u}\cdot\nabla{f}=0$
$$
\omega^{ik}Q^{kj}\vec{v}^i\vec{v}^j=\omega^{ik}\big(Q^{kj}\vec{v}^j\big)\vec{v}^i=\lambda_1\omega^{ik}\vec{v}^k\vec{v}^{i}=0,
$$
$$
Q^{ik}\omega^{kj}\vec{v}^i\vec{v}^j=\big(Q^{ik}\vec{v}^i\big)\omega^{kj}\vec{v}^j=\lambda_1\omega^{kj}\vec{v}^k\vec{v}^j=0
$$
 hence we get
\begin{align}
\partial_tf&=\Delta{f}-a\lambda_1-c\tr(Q^2)\lambda_1+b\Big[\lambda_1^2-\dfrac{\tr(Q^2)}{3}\Big]\nonumber\\
&=\Delta{f}-\lambda_1\big[a+c(\lambda_1^2+\lambda_2^2+\lambda_3^2)\big]+b\Big(\lambda_1^2-\dfrac{\lambda_1^2+\lambda_2^2+\lambda_3^2}{3}\Big)
\quad\mbox{at } (t_0, x_0).\label{scalar equ}
\end{align}
Using Cauchy-Schwarz inequality and $Q\in\sS_0^{(3)}$, we get
$$
  \lambda_1^2+\lambda_2^2+\lambda_3^2\geq\lambda_1^2+\dfrac{(\lambda_2+\lambda_3)^2}{2}=\dfrac{3}{2}\lambda_1^2,
$$
which combined with \eqref{scalar equ} at $(t_0, x_0)$ gives
\begin{align*}
\partial_tf|_{(t_0,x_0)}&\leq\Delta{f}-a\lambda_1-\dfrac{3c}{2}\lambda_1^3+\dfrac{b}{2}\lambda_1^2\Big|_{(t_0,x_0)}
\leq-\dfrac{3c}{2}\lambda_1\Big(\lambda_1^2-\dfrac{b}{3c}\lambda_1+\dfrac{2a}{3c}\Big)\Big|_{(t_0,x_0)}\\
&<0.
\end{align*}
Above in the last inequality we used \eqref{large assumption}.
However, \eqref{maximum f} indicates that
$$
  \partial_tf|_{(t_0,x_0)}\geq 0,
$$
which is a contradiction. Therefore, we conclude that
\begin{equation}\label{upper bound}
\lambda_1(t, x)\leq \frac{b+\sqrt{b^2-24ac}}{6c}, \quad\forall (t,
x)\in (0,T]\times\Omega.
\end{equation}

\noindent{\bf Step 2.} Let
\begin{equation}\label{minimum point}
\lambda_3(\tilde{t}, \tilde{x})=\displaystyle\min_{(t, x)\in[0,
T]\times{\bar\Omega}}\lambda_3(t, x).
\end{equation}
We shall again show that $ \lambda_3(\tilde{t},\tilde{x}) \geq
\dfrac{-b-\sqrt{b^2-24ac}}{12c}, $ by contradiction. Suppose
\begin{equation}\label{small assumption}
(\tilde{t}, \tilde{x})\in(0, T]\times\Omega,
\;\mbox{and}\;\lambda_3(\tilde{t},\tilde{x})<\frac{-b-\sqrt{b^2-24ac}}{12c}.
\end{equation}
Let $\vec{w}\in\mathbb{S}^2$ be the corresponding unit eigenvector,
such that
$Q(\tilde{t},\tilde{x})\vec{w}=\lambda_3(\tilde{t},\tilde{x})\vec{w}$.
Meanwhile, we denote
$$
  g(t, x)=\langle Q(t, x)\vec{w}, \vec{w}\rangle_{\mathbb{R}^3},
$$
then we see from \eqref{initial setting} and \eqref{minimum point}
that
\begin{equation}\label{minimum g}
g(\tilde{t},\tilde{x})=\displaystyle\min_{(t, x)\in[0,
T]\times{\bar\Omega}}g(t, x)
\end{equation}
After taking the matrix inner product of equation \eqref{Q-equ} with
$\vec{w}\vec{w}^t$, and evaluating at $(\tilde{t},\tilde{x})$, it
gives
\begin{align*}
\mathbf{u}\cdot\nabla{g}&=0,\\
\omega^{ik}Q^{kj}\vec{w}^i\vec{w}^j=\omega^{ik}\big(Q^{kj}\vec{w}^j\big)\vec{w}^i=\lambda_1\omega^{ik}\vec{w}^k\vec{w}^j&=0,\\
Q^{ik}\omega^{kj}\vec{w}^i\vec{w}^j=\big(Q^{ik}\vec{w}^i\big)\omega^{kj}\vec{w}^j=\lambda_1\omega^{kj}\vec{w}^k\vec{w}^j&=0.
\end{align*}
Consequently, we obtain
\begin{align}
\partial_tg&=\Delta{g}-a\lambda_3-c\tr(Q^2)\lambda_3+b\Big[\lambda_3^2-\dfrac{\tr(Q^2)}{3}\Big]\nonumber\\
&=\Delta{g}-\lambda_3\big[a+c(\lambda_1^2+\lambda_2^2+\lambda_3^2)\big]+b\Big(\lambda_3^2-\dfrac{\lambda_1^2+\lambda_2^2+\lambda_3^2}{3}\Big)
\quad\mbox{at } (\tilde{t},\tilde{x}).\label{scalar equ 2}
\end{align}
We claim
\begin{equation}\label{eqn:claim}
\partial_tg\big|_{(\tilde{t},\tilde{x})}
\geq-\dfrac{3c}{2}\lambda_3\Big(\lambda_3^2+\dfrac{b}{6c}\lambda_3+\dfrac{a}{6c}\Big)\Big|_{(\tilde{t},\tilde{x})}.
\end{equation}
Here, we focus on the proof of the theorem and the proof of the
claim will be postponed to the end of the section. Combining
\eqref{small assumption}, and the claim \eqref{eqn:claim}, we get
$$
  \partial_tg|_{(\tilde{t},\tilde{x})}>0.
$$
On the other hand, however, based on \eqref{minimum g} one can
deduce that
$$
  \partial_tg|_{(\tilde{t},\tilde{x})}\leq 0,
$$
which is again a contradiction. Thus
\begin{equation}\label{lower bound}
\lambda_3(t, x)\geq -\frac{b+\sqrt{b^2-24ac}}{12c}, \quad\forall (t,
x)\in (0,T]\times\Omega.
\end{equation}
The proof is complete by combining \eqref{upper bound} and
\eqref{lower bound}.
\end{proof}

It remains to prove the claim \eqref{eqn:claim}.
\begin{proof}
We divide the proof of \eqref{eqn:claim} into three cases. By \eqref{eigsumzero}
and \eqref{minimum point}, we know that
$|\lambda_2|\leq |\lambda_3|$.

\smallskip

\noindent\textbf{Case 1:} $\lambda_2|_{(\tilde{t},\tilde{x})}\geq 0$

\smallskip

First of all, note that
$$
  \dfrac{3}{2}\lambda_3^2=\frac{(\lambda_1+\lambda_2)^2}{2}+\lambda_3^2 \leq\lambda_1^2+\lambda_2^2+\lambda_3^2\leq(\lambda_1+\lambda_2)^2+\lambda_3^2=2\lambda_3^2,
$$
which together with \eqref{scalar equ 2} and the fact that $g$
attains minimum at $(\tilde{t},\tilde{x})$ yields
\begin{eqnarray*}
\partial_tg|_{(\tilde{t},\tilde{x})}
\geq\Delta{g}-a\lambda_3-\dfrac{3c}{2}\lambda_3^3+\frac{b}{3}\lambda_3^2\Big|_{(\tilde{t},\tilde{x})}
&\geq&-a\lambda_3-\dfrac{3c}{2}\lambda_3^3+\frac{b}{3}\lambda_3^2\Big|_{(\tilde{t},\tilde{x})}\\
&\geq&-\dfrac{3c}{2}\lambda_3\Big(\lambda_3^2+\dfrac{b}{6c}\lambda_3+\dfrac{a}{6c}\Big)\Big|_{(\tilde{t},\tilde{x})}
\end{eqnarray*}

\noindent\textbf{Case 2:}
$\dfrac{(\sqrt{5}-1)}{2}\lambda_3\big|_{(\tilde{t},\tilde{x})}\leq\lambda_2|_{(\tilde{t},\tilde{x})}<0$

\smallskip

In this case, again thanks to \eqref{eigsumzero} $|\lambda_1|>|\lambda_3|$ at ${(\tilde{t},\tilde{x})}$
we have
\begin{eqnarray*}
  2\lambda_3^2<\lambda_1^2+\lambda_2^2+\lambda_3^2&=&(\lambda_2+\lambda_3)^2+\lambda_2^2+\lambda_3^2=2(\lambda_2^2+\lambda_3^2+\lambda_2\lambda_3)\\
  &=&2(\lambda_2^2+\lambda_2\lambda_3-\lambda_3^2)+4\lambda_3^2\leq
  4\lambda_3^2.
\end{eqnarray*}
Hence
\begin{eqnarray*}
\partial_tg|_{(\tilde{t},\tilde{x})}\geq\Delta{g}-a\lambda_3-2c\lambda_3^3-\frac{b}{3}\lambda_3^2\Big|_{(\tilde{t},\tilde{x})}
&=&-2c\lambda_3\Big(\lambda_3^2+\dfrac{b}{6c}\lambda_3+\dfrac{a}{2c}\Big)\Big|_{(\tilde{t},\tilde{x})}\\
&\geq&-\dfrac{3c}{2}\lambda_3\Big(\lambda_3^2+\dfrac{b}{6c}\lambda_3+\dfrac{a}{6c}\Big)\Big|_{(\tilde{t},\tilde{x})}
\end{eqnarray*}

\noindent\textbf{Case 3:}
$\lambda_3|_{(\tilde{t},\tilde{x})}\leq\lambda_2|_{(\tilde{t},\tilde{x})}<\dfrac{(\sqrt{5}-1)}{2}\lambda_3\big|_{(\tilde{t},\tilde{x})}
< 0$

\smallskip

Note that
\begin{eqnarray*}
   4\lambda_3^2&<&2(\lambda_2^2+\lambda_2\lambda_3-\lambda_3^2)+4\lambda_3^2=2(\lambda_2^2+\lambda_3^2+\lambda_2\lambda_3)\\
 & =&(\lambda_2+\lambda_3)^2+\lambda_2^2+\lambda_3^2
=\lambda_1^2+\lambda_2^2+\lambda_3^2 \quad\mbox{at } (\tilde{t},\tilde{x}),
\end{eqnarray*}
which gives
\begin{align*}
&-\lambda_3\big[a+c(\lambda_1^2+\lambda_2^2+\lambda_3^2)\big]+b\Big(\lambda_3^2-\dfrac{\lambda_1^2+\lambda_2^2+\lambda_3^2}{3}\Big)\\
&=-2c\lambda_3\left[(\lambda_2^2+\lambda_3^2+\lambda_2\lambda_3)+\frac{b}{6c}\frac{(2\lambda_2^2+2\lambda_2\lambda_3-\lambda_3^2)}{\lambda_3}
+\frac{a}{2c}\right]\\
&\defeq -2c\lambda_3H(\lambda_2,\lambda_3)  \quad\mbox{at }
(\tilde{t},\tilde{x}).
\end{align*}
We proceed to show that
\begin{equation}\label{equivalent claim}
  H(\lambda_2,\lambda_3)\geq \lambda_3^2+\dfrac{b}{6c}\lambda_3+\dfrac{a}{6c} \quad\mbox{at }
(\tilde{t},\tilde{x}),
\end{equation}
which is equivalent to
\begin{equation}\label{equivalent claim-1}
\lambda_2^2+\lambda_2\lambda_3+\frac{b}{3c}\frac{(\lambda_2^2+\lambda_2\lambda_3-\lambda_3^2)}{\lambda_3}
+\frac{a}{3c}\geq 0  \quad\mbox{at } (\tilde{t},\tilde{x}).
\end{equation}
Let us denote
$\mu=\lambda_2^2+\lambda_2\lambda_3-\lambda_3^2\big|_{(\tilde{t},\tilde{x})}$,
then $0<\mu\leq\lambda_3^2$ due to the assumption in {\bf Case 3}
and \eqref{equivalent claim-1}
 is reduced to
\begin{equation}\label{equivalent claim-2}
\Big(1+\dfrac{b}{3c}\dfrac{1}{\lambda_3}\Big)\mu+\lambda_3^2+\frac{a}{3c}\geq
0 \quad\mbox{at } (\tilde{t},\tilde{x}).
\end{equation}
By \eqref{small assumption}, we have
\begin{align}
-3\leq
1-\dfrac{4b}{b+\sqrt{b^2-24ac}}\leq\Big(1+\dfrac{b}{3c}\dfrac{1}{\lambda_3}\Big)_{(\tilde{t},\tilde{x})}\leq
1
\end{align}
If
$0\leq\Big(1+\dfrac{b}{3c}\dfrac{1}{\lambda_3}\Big)_{(\tilde{t},\tilde{x})}\leq
1$, then \eqref{equivalent claim-2} is automatically true.
 Otherwise
since $\mu$ is a monotone decreasing, nonnegative function of
$\lambda_2$ on the given interval, we have
\begin{eqnarray*}
\Big(1+\dfrac{b}{3c}\dfrac{1}{\lambda_3}\Big)\mu+\lambda_3^2+\frac{a}{3c}
&\geq&\Big(1+\dfrac{b}{3c}\dfrac{1}{\lambda_3}\Big)\lambda_3^2+\lambda_3^2+\frac{a}{3c}\\
&=&2\Big(\lambda_3^2+\dfrac{b}{6c}\lambda_3+\dfrac{a}{6c}\Big)>0
\quad\mbox{at } (\tilde{t},\tilde{x}),
\end{eqnarray*}
where we used \eqref{small assumption} in the last inequality above.
In all, \eqref{equivalent claim-2} is valid, and so is
\eqref{equivalent claim}. As a consequence,
\begin{align*}
\partial_tg|_{(\tilde{t},\tilde{x})}\geq -2c\lambda_3\Big(\lambda_3^2+\dfrac{b}{6c}\lambda_3+\dfrac{a}{6c}\Big)\Big|_{(\tilde{t},\tilde{x})}
\geq-\frac{3c}{2}\lambda_3\Big(\lambda_3^2+\dfrac{b}{6c}\lambda_3+\dfrac{a}{6c}\Big)\Big|_{(\tilde{t},\tilde{x})}.
\end{align*}
The proof of the claim \ref{eqn:claim} is complete.
\end{proof}
\bigskip
{After this, Corollary \ref{Corollary} can be easily established.}

\medskip

{\it Proof of Corollary \ref{Corollary}. }
Without loss of generality we may assume
$$
Q(t, x)\in C^{1,2}((0,T)\times\Omega)\cap C([0, T]\times\bar\Omega),
$$
and let
\begin{equation}
\lambda_1(t_0, x_0)=\displaystyle\max_{(t, x)\in[0,
T]\times{\bar\Omega}}\lambda_1(t, x).
\end{equation}
If
$$
  \lambda_1(t_0,x_0)>\max\Big[\frac{b+\sqrt{b^2-24ac}}{6c},
\displaystyle\max_{\bar{\Omega}}\lambda_1(Q_0) \Big],
$$
then $t_0\in (0, T]$ and $x_0\in\Omega$. As a consequence, it
follows from the same argument as in the proof of Theorem \ref{main
theorem} that
$$
  0\leq
  -\dfrac{3c}{2}\lambda_1\Big(\lambda_1^2-\dfrac{b}{3c}\lambda_1+\dfrac{2a}{3c}\Big)\Big|_{(t_0,x_0)}<0
$$
due to the assumption that $
\lambda_1(t_0,x_0)>\frac{b+\sqrt{b^2-24ac}}{6c}$, which is a
contradiction. The corresponding lower bound for $\lambda_3(t, x)$
can be proved in a similar manner.

\qed

\begin{remark}
With minor modifications, one may check easily that the above
arguments are also valid for the whole space case that is shown in
\cite{XZ16}.
\end{remark}

\section{Appendix}

In this appendix section using the same idea as \cite{XZ16}, we
prove

\begin{lemma}\label{lemma on regularizing effect}

Let $\mathbf{u}\in H^1(0, T; H^1(\Omega))\cap L^\infty(0, T;
H^2(\Omega))$, $\mathbf{u}_\delta\in C^\infty
([0,T]\times\bar{\Omega})$,
$\nabla\cdot\mathbf{u}=\nabla\cdot\mathbf{u}_\delta=\mathbf{u}|_{\partial\Omega}=\mathbf{u}_\delta|_{\partial\Omega}=0$
be such that $\mathbf{u}_\delta \rightarrow \mathbf{u}$ as
$\delta\rightarrow 0$ strongly in $H^1(0, T; H^1(\Omega))\cap
L^\infty(0, T; H^2(\Omega))$. Let $\dQ$ be the unique classical
solution in $C^{1,2}((0,T)\times\Omega)\cap C([0,
T]\times\bar\Omega)$ of the system
\begin{eqnarray}
&&\dQ_t+\mathbf{u_\delta}\cdot\nabla\dQ-\omega_\delta\dQ+\dQ\omega_\delta\\
&&=\Gamma\left(L\Delta\dQ-a\dQ+b\Big[(\dQ)^2-\frac{\tr(\dQ)^2}{3}\mathbb{I}\Big]-c\dQ\tr(\dQ)^2\right),\non
\end{eqnarray}
\begin{eqnarray}
\dQ(t, x)|_{\partial\Omega}=\tilde{Q}(x),
\end{eqnarray}
where
$\omega_\delta=\frac{\nabla\mathbf{u}^\delta-\nabla^T\mathbf{u}^\delta}{2}$.
Assume that
\begin{equation}
\tilde{m}\leq\lambda_i (\dQ(t,x))\leq\tilde{M}, \qquad \forall\,
1\leq i\leq 3,\quad (t, x)\in [0, T]\times\Omega.
\end{equation}
Then $Q^{(\delta)}(t,x)\rightarrow Q(t,x)$ as $\delta \rightarrow
0$, $\forall t\ge 0, x\in\Omega$, where $Q$ is the unique solution
in $H^1(0, T; H^2(\Omega))\cap L^\infty(0, T; H^3(\Omega))$ of
\begin{align}
Q_t+\mathbf{u}\cdot\nabla{Q}-\omega{Q}+Q\omega&=\Gamma\Big(L\Delta{Q}-aQ+b\Big[Q^2-\frac{\tr(Q^2)}{3}\mathbb{I}\Big]-cQ\tr(Q^2)\Big),\label{Q-equ+}\\
Q(t, x)|_{\partial\Omega}&=\tilde{Q}(x).
\end{align}
Furthermore, we have an eigenvalue constraint on $Q$ such that
\begin{equation}
\tilde{m}\le \lambda_i (Q(t,x))\leq\tilde{M},\qquad \forall\, 1\leq
i\leq 3,\quad (t, x)\in (0, T]\times\Omega \label{eigenvalue-const}
\end{equation}
provided the initial data $Q_0$ and boundary data $\tilde{Q}$ have
the same constraint.
\end{lemma}
\begin{remark}
{ Theorem 1.1 in \cite{LY16} is only one
step away from the existence of local classical solutions to the
evolution problem \eqref{NSE}-\eqref{IC-BC}, which nevertheless
cannot be improved using the method therein. However, it remains to
be an interesting question to study.}
\end{remark}

\begin{proof}
Without loss of generality, we set $\Gamma=L=1$. To begin with, we
show apriori $L^\infty$ bound on $Q$. Since $c>0$, there exists
$\eta_0>0$, such that
$$
 -a|M|^2+b\tr(M)^3-c|M|^4\le 0, \qquad\forall\, |M|\geq\eta_0,
 \,M\in\mathcal{S}_0^3
$$
Let $\eta\defeq\max\{\|Q_0\|_{L^\infty(\Omega)},
\|\tilde{Q}\|_{L^\infty(\Omega)}, \eta_0\}$. Multiplying
\eqref{Q-equ+} with $Q(|Q|^2-\eta)_+$, then integrating over
$\Omega$ and using integration by parts, we obtain
\begin{align*}
\frac{1}{2}\frac{d}{dt}\int_{\Omega} (|Q|^2-\eta)_+^2
&=-\int_{\Omega}|\nabla Q|^2(|Q|^2-\eta)_+\,dx-\int_{\Omega} |\nabla (|Q|^2-\eta)_+|^2\,dx\\
&\qquad+\int_{\Omega} (-a|Q|^2+b\tr(Q^3)-c|Q|^4)(|Q|^2-\eta)_+\le 0
\end{align*}
Thus
\begin{equation}\label{Qinftyeta}
\|Q(t,\cdot)\|_{L^\infty(\Omega)}\le \eta,\quad \forall\, 0\leq
t\leq T.
\end{equation}
In the same way, we conclude
\begin{equation}\label{dQinftyeta}
\|\dQ(t,\cdot)\|_{L^\infty(\RR^3)}\le \eta, \quad\forall\, 0\leq
t\leq T,\, \forall\delta>0.
\end{equation}
Let $\dR\defeq\dQ-Q\in\mathcal{S}_0^3$. Then it is easy to see
\begin{eqnarray*}
&&\partial_t \dR+\mathbf{u}_\delta\nabla R_\delta-\omega_\delta \dR+\dR\omega_\delta\\
&&=\Delta\dR-a\dR+(\mathbf{u}-\mathbf{u}_\delta)\nabla Q-(\omega-\omega_\delta)Q+Q(\omega-\omega_\delta)\\
&&+b\Big[\dR\dQ+Q\dR-\frac{\tr(\dR\dQ+Q\dR)}{3}\mathbb{I}\Big]\\
&&-c\Big[|\dQ|^2\dR+\tr(\dR\dQ+Q\dR)Q\Big],
\end{eqnarray*}
with initial and boundary datum $\dR_0=\dR|_{\partial\Omega}\equiv
0$. Multiplying the above equation with $\dR$, integrating over
$\Omega$, by \eqref{Qinftyeta} and \eqref{dQinftyeta}, we get
\begin{eqnarray}
&&\frac{d}{dt}\int_{\Omega} |\dR|^2\,dx\non\\
&&\leq C\int_{\Omega}|\dR|^2\,dx+C\int_{\Omega}
\big[(\mathbf{u}-\mathbf{u}_\delta)\nabla Q-(\omega-\omega_\delta)Q+Q(\omega-\omega_\delta)\big]\dR\,dx\non\\
&&\leq C\|\dR\|_{L^2}^2+C\big(\|\nabla Q\|_{L^2}
\|\mathbf{u}-\mathbf{u}_\delta\|_{L^2}+\|Q\|_{L^2}\|\omega(s)-\omega_\delta(s)\|_{L^2}\big)\non\\
&&\leq
C\|\dR\|_{L^2}^2+C\big(\|\mathbf{u}-\mathbf{u}_\delta\|_{L^2}+\|\omega(s)-\omega_\delta(s)\|_{L^2}\big).\non\\
\end{eqnarray}
Hence Gronwall's inequality gives
\begin{align*}
  \|\dR(t,\cdot)\|_{L^2}^2\leq
  Ce^{Ct}\int_0^t\|\mathbf{u}-\mathbf{u}_\delta\|_{H^1}(s)\,ds\rightarrow
  0, \quad\mbox{as } \delta\rightarrow 0
\end{align*}
due to the assumption that $\mathbf{u}_\delta \rightarrow
\mathbf{u}$ strongly in $L^\infty(0, T; H^2(\Omega))$. Therefore,
combined with the fact that $\dQ, Q\in C([0, T]\times\bar\Omega)$ we
get
\begin{equation}\label{pointwise convergence}
  Q^\delta(t, x)\rightarrow Q(t, x), \quad\forall\, (t, x)\in [0, T]\times\Omega
\end{equation}
Moreover, \eqref{eigenvalue-const} follows from \eqref{pointwise
convergence} and Lemma \ref{continuous eigenvalue}

\end{proof}

\bigskip

\noindent \textbf{Acknowledgements.} We thank the anonymous referees
for their careful reading and useful suggestions to improve our
paper, {especially the observation that leads to Corollary
\ref{Corollary}, which can be considered an improved result of
Theorem \ref{main theorem}}. The work of A. Contreras was partially
supported by a grant from the Simons Foundation \# 426318. And the
work of W.~J. Zhang is supported by the start up fund from
Department of Mathematics, Rutgers university. We want to thank our
friends Xavier Lamy, Yuning Liu and Arghir Zarnescu for their kind
discussions. In particular, X. Xu would like to express his
gratitude to Arghir for his consistent support and academic
communications over a series of topics on the mathematical
$Q$-tensor theory during the past six years. Without his idea
proposed in \cite{XZ16}, this paper would not come out.


\end{document}